\documentclass[11pt]{article}
\usepackage{epsfig}
\usepackage{t1enc}
    \usepackage[latin1]{inputenc}
    \usepackage[english]{babel}
        \usepackage{latexsym}
\usepackage{amssymb}
\usepackage{amsmath}
\usepackage{euscript}
\newcommand{\bz}{\mathbb{Z}}
\newcommand{\br}{\mathbb{R}}
\newcommand{\bs}{\textbf{S}}

\newcommand{\bo}{\textbf{0}}
\newcommand{\bb}{\textbf{B}}

\def\begg{\begin{equation}}
\def\endd{\end{equation}}

\newcommand{\bp}{{\bf P}}

\begin{document}
\setlength{\arraycolsep}{.136889em}
\renewcommand{\theequation}{\thesection.\arabic{equation}}
\newtheorem{thm}{Theorem}[section]
\newtheorem{propo}{Proposition}[section]
\newtheorem{lemma}{Lemma}[section]
\newtheorem{corollary}{Corollary}[section]
\newtheorem{remark}{Remark}[section]

\smallskip

\centerline{\Large\bf LIMIT THEOREMS FOR LOCAL AND OCCUPATION
TIMES}
\smallskip
\centerline{\Large\bf OF RANDOM WALKS AND BROWNIAN MOTION ON A SPIDER}

\bigskip\bigskip
\bigskip\bigskip
\renewcommand{\thefootnote}{1}

 \bigskip\bigskip
 \renewcommand{\thefootnote}{1}
 \noindent {\textbf{Endre Cs\'{a}ki}\footnote{Research supported by
Hungarian National Research, Development and Innovation Office-NKFIH Grant
 No. K 108615.}}
 \newline
 Alfr\'ed R\'enyi Institute of Mathematics, Hungarian Academy of
 Sciences, Budapest, P.O.B. 127, H-1364, Hungary. E-mail address:
 csaki.endre@renyi.mta.hu

 \bigskip
 \renewcommand{\thefootnote}{2}
 \noindent {\textbf{Mikl\'os Cs\"org\H{o}}\footnote{Research
 supported by an NSERC Canada Discovery Grant at Carleton
 University}}
 \newline
 School of Mathematics and Statistics, Carleton University, 1125
 Colonel By Drive, Ottawa, Ontario, Canada K1S 5B6. E-mail address:
 mcsorgo@math.carleton.ca

\bigskip
\renewcommand{\thefootnote}{3}
\noindent
{\textbf{Ant\'{o}nia
F\"{o}ldes}\footnote{Research supported by a PSC CUNY Grant, No.
69040-0047.}}
\newline
Department of Mathematics, College of Staten
Island, CUNY, 2800 Victory Blvd., Staten Island, New York 10314,
U.S.A.  E-mail address: antonia.foldes@csi.cuny.edu

\bigskip
\renewcommand{\thefootnote}{4}
\noindent
{\textbf{P\'al R\'ev\'esz}\footnote{Research supported by
Hungarian National Research, Development and Innovation Office-NKFIH Grant
 No. K 108615.}}
\newline
Institut f\"ur
Statistik und Wahrscheinlichkeitstheorie, Technische Universit\"at
Wien, Wiedner Hauptstrasse 8-10/107 A-1040 Vienna, Austria.
E-mail address: revesz.paul@renyi.mta.hu
\bigskip

\medskip
\noindent{\bf Abstract}\newline
A simple random walk and a Brownian motion are considered on a spider that is 
a collection of half lines (we call them legs) joined at the origin. We give 
a strong approximation of these two objects and their local times. For fixed 
number of legs we establish limit theorems for $n$ step local and occupation 
times.


\medskip

\vspace{.1cm}

\medskip

\section{Introduction }
\renewcommand{\thesection}{\arabic{section}} \setcounter{equation}{0}
\setcounter{thm}{0} \setcounter{lemma}{0}

Consider the following collection of half lines on the complex plane 
$$\textbf{SP}(N):=\{v_N(x,j),\, x\geq 0,\, j=1,2,\ldots,N\},$$
where, with $i=\sqrt{-1}$, 
$$
v_N(x,j)=x\exp\left(\frac{2\pi i(j-1)}{N}\right).
$$
We will call $\textbf{SP}(N)$ a spider with $N$ legs. Also,

$$v_N(0):=v_N(0,1)=v_N(0,2)=...=v_N(0,N)$$
\noindent
is called the body of the spider, and
$L_j:=\{v_N(x,j),\, x> 0 \}$ is the $j$-th leg of the spider.

We define the distance on $\textbf{SP}(N)$ by
$$
|v_N(x,j)-v_N(y,j)|=|x-y|, \quad j=1,\ldots,N
$$
$$
|v_N(x,j)-v_N(y,k)|=x+y, \quad j,k=1,\ldots,N,\, \, j\neq k.
$$

In this paper the number of legs $N$ is fixed, so we will often suppress 
$N$ in the notation and, instead of $v_N(x,j)$ or $v_N(0)$, we will simply 
write $v(x,j)$ or $v(0)=\bo$, whenever convenient.

We consider a random walk $\bs_n,\,\, n=0,1,2 \ldots$, on
$\textbf{SP}(N)$ that starts from the body of the spider, i.e.,
$\bs_0=v_N(0)=\bo,$ with the following transition probabilities:
$${\bf P} (\bs_{n+1}=v_N(1,j)|\bs_n=v_N(0))=p_j,\quad j=1,...,N, $$
with
$$\sum_{j=1}^N p_j=1,$$
and, for $r=1,2,...,\quad j=1,...,N$,
$${\bf P} (\bs_{n+1}=v_N(r+1,j)|\bs_n=v_N(r,j))={\bf P}
(\bs_{n+1}=v_N(r-1,j)|\bs_n=v_N(r,j))=\frac{1}{2}.$$
The thus defined random walk $\bs_n$ on $\textbf{SP}(N)$ will be called 
random walk on spider (RWS) or simply spider walk in this paper.

Observe that the particular case $N=2$, $p_1=p_2=1/2$ corresponds to the
simple symmetric random walk $S(n),\, n=0,1,2,\ldots,$ on the line. 
The spider walk $\bs_n$ can be constructed from $S(n)$, $n=0,1,\ldots$, 
as follows. Consider the absolute value $|S(n)|,\, n=1,2,\ldots$, that 
consists of infinitely many excursions from zero, and let $G_1,G_2,\ldots$, 
denote the excursion intervals. Put these excursions, independently of each 
other, on leg $j$ of the spider with probability $p_j$, $j=1,2,\ldots,N$. 
What we obtain this way is the first $n$ steps of the spider walk 
$\bf{S}_\cdot,$ created from  the first $n$ steps of the random walk 
$S(\cdot)$. Let $\gamma_i$, $i=1,2,\ldots$, be i.i.d. random variables with
$$
P(\gamma_i=j)=p_j, \quad j=1,2,\ldots,N,
$$
that are independent from the simple symmetric random walk $S(\cdot)$. 
Consequently, one can also  define the just introduced spider walk as 
${\bf S}_0=\bf{0},$
\begin{equation}
\bs_n:=\sum_{m=1}^\infty I\{n\in G_m\}v_N(|S(n)|,\gamma_m),\, \, n=1,2,\ldots,
\label{RWSdef}
\end{equation}
if $S(n)\neq 0$, and $\bs_n=v_N(0)$ if $S(n)=0$. 

In view of this definition, and the notations already used, in the sequel
$\bs_n$ will stand for a spider walk and $S(n)$ for a simple symmetric random 
walk on the line, with respective probabilities denoted by ${\bf P}$ and $P$.

The limit process is the so-called Brownian motion on spider (BMS), or 
simply Brownian spider, a  version of Walsh's Brownian motion (cf Walsh 
\cite{WAL}, and Introduction of Cs\'aki et al. \cite{Cs14}), that can be 
constructed from a standard Brownian motion 
$\{B(t),\, t\geq 0\}$ on the line as follows. The process 
$\{|B(t)|\, t\geq 0\}$ has a countable number of excursions from zero, and 
denote by $J_1,J_2,\ldots$, a fixed enumeration of its excursion intervals 
away from zero. Then, for any $t>0$ for which $B(t)\neq 0$, we have that 
$t\in J_m$ for one of the values of $m=1,2,\ldots$. Let
$\kappa_m, \, m=1,2,\ldots,$ be i.i.d. random variables, independent of 
$B(\cdot)$ with
$$
P(\kappa_m=j)=p_j,\quad j=1,2,\ldots, N.
$$
We now construct the Brownian spider $\{{\bb}(t),\, t\geq 0\}$ by putting 
the excursion whose interval is $J_m$ to the $\kappa_m$-th leg of the spider 
$\textbf{SP}(N)$. 
Hence we can define the Brownian spider ${\bb}(\cdot)$ by
\begg
{\bb}(t):=\sum_{m=1}^\infty I\{t\in J_m\}v_N(|B(t)|,\kappa_m), \quad
{\rm if}\quad B(t)\neq 0, \label{sdef}
\endd
and
$$
{\bb}(t):=v_N(0)=\bo,\quad {\rm if}\quad B(t)=0.
$$

In \cite{Cs14} we investigated the heights of the spider walk ${\bf S}_n$
on the legs of the spider. We also established a strong approximation as 
follows. 
\begin{thm} On a rich enough probability space one can define a {\rm BMS}
$\{{\bf B}(t),\, t\geq 0\}$ and an {\rm RWS} 
$\{{\bf S}_n,\, n=0,1,2,\ldots\}$, both on ${\bf SP}(N)$, and both selecting 
their legs  with the same probabilities $p_j,\, j=1,2,\ldots,N$, so that, 
as $n\to\infty$, we have
$$
|{\bf S}_n-{\bf B}(n)|=O(n^{1/4}(\log n)^{1/2}(\log\log n)^{1/4}) \quad a.s.
$$
\end{thm}

In this paper we keep $N$ and $p_1,\ldots,p_N$ fixed and consider limit 
theorems concerning local and occupation times on the legs, as the number 
of steps $n$ tends to infinity. 

\section{Local times}
\renewcommand{\thesection}{\arabic{section}} \setcounter{equation}{0}
\setcounter{thm}{0} \setcounter{lemma}{0}

Throughout we use the notation $I\{A\}$ for the indicator of an event $A$
in the brackets, i.e., $I\{A\}=1$ if $A$ occurs and $I\{A\}=0$ otherwise.
  
The local time of an RWS $\bs_n$ on ${\bf SP}(N)$ is defined as 
$$
\boldsymbol{\xi}((r,j),n):=\sum_{i=1}^n I\{\bs_i=v(r,j)\},\, \, r=1,2,\ldots, 
\, j=1,2,\ldots,N,
$$
$$
\boldsymbol{\xi}(0,n):=\sum_{i=1}^n I\{{\bs}_i=v(0)\},
$$
$n=1,2,\ldots$.

The local time of a BMS is defined by
$$
\boldsymbol{\eta}((x,j),t):=\lim_{\varepsilon\to 0}\frac{1}{2\varepsilon}
\int_0^t I\{{\bb}(s)\in (v(x-\varepsilon,j),v(x+\varepsilon,j))\}\, ds, 
\, \,x>0,\qquad j=1,2,\ldots N,
$$
$$
\boldsymbol{\eta}(0,t):=\lim_{\varepsilon\to 0}\frac{1}{2\varepsilon}
\int_0^t \sum_{j=1}^N I\{{\bb}(s)\in (0,v(\varepsilon,j))\}\, ds.
$$

Note that by the constructions given in Section 1 we have
$$
\boldsymbol{\xi}(0,n)=\xi(0,n)\quad {\rm and} \quad
\boldsymbol{\eta}(0,n)=\eta(0,n),
$$
where $\xi(0,n)$ and $\eta(0,n)$, are the local times at zero of the simple 
symmetric random walk $S(\cdot)$ and standard Brownian motion $B(\cdot)$,
respectively, that are used in the constructions. 

First we consider the case $N=2$, i.e., the so called skew random walk 
(SRW) and skew Brownian motion (SBM). In this case we let 
$p_1=p$, $p_2=1-p=q$.

\subsection{Skew random walk and skew Brownian motion}

The skew random walk (SRW) $\{{\bf S}_n^*,\, n=0,1,2,\ldots\}$ with parameter 
$p$ is the particular case of an RWS with $N=2$. It is a 
Markov chain on ${\bz}$ with transition probabilities
$$
{\bf P}(0,1)=p, \,\, {\bf P}(0,-1)=q=1-p,\,\, 
{\bf P}(x,x+1)={\bf P}(x,x-1)=1/2,\,\, x=\pm 1,\pm 2,\ldots
$$
In this case the construction  yielding (\ref{RWSdef}) is equivalent to
${\bf S}_n^*:=0 ,$ \, if\, $S(n)=0,$ and
\begin{equation}
{\bf S}_n^*:=\sum_{m=1}^\infty I\{n\in G_m\}\delta_m|S(n)|, \qquad {\rm if}  
\quad S(n)\neq 0, \,\,n=1,2,\dots,
\label{SRWdef}
\end{equation}
where $\delta_m, \, \, m=1,2,\ldots$ are i.i.d. random variables with
$P(\delta_1=1)=p=1-P(\delta_1=-1)$, and $G_1,G_2,\ldots$ are excursion 
intervals of $S(\cdot)$. 

The skew Brownian motion (SBM) $\bb^*(\cdot)$ with parameter $p$ can be 
defined as follows (cf. \cite{ABT}). Let $B(\cdot)$ be a standard Brownian 
motion on the line, and let $J_1,J_2,\ldots$ be the excursion intervals from 
zero of $|B(\cdot)|$. Let $\delta_i,\, \, i=1,2,\ldots$, be i.i.d. random 
variables, independent of $B(\cdot)$, with $P(\delta_1=1)=p=1-P(\delta_1=-1)$.
Then
\begin{equation}
{\bb^*}(t):=\sum_{m=1}^\infty I\{t\in J_m\}\delta_m|B(t)|, \quad
{\rm if}\quad B(t)\neq 0, 
\label{sdef2}
\end{equation}
\begin{equation}
{\bb^*}(t):=0,\qquad {\rm if}\quad B(t)=0.
\end{equation}

The skew Brownian motion was introduced by It\^o and McKean \cite{IM} 
(cf. also the Introduction of Harrison and Shepp  \cite{HS}).
Weak invariance principle between skew random walk and skew Brownian motion was 
established by Harrison and Shepp \cite{HS} and Cherny et al. \cite{CSY}.
We note that throughout this paper we use the above introduced respective 
abbreviations SRW and SBM. For further results concerning SRW and SBM, we 
refer to Revuz and Yor \cite{RY}, Lejay \cite{Le}, \cite{Le2}, Appuhamillage 
et al. \cite{ABT}, Hajri \cite{Ha} and references given in these papers.     

We note in passing that the above definition of SBM is a special case of 
BMS as in (\ref{sdef}) with $N=2$, $p_1=p$, $p_2=1-p=q$, so that, 
in this case, the first one of the two legs is the positive half-line, and the 
second one is the negative half-line. In particular, we thus have 
$$
v_2(|B(t)|,\kappa_m)=\delta_m|B(t)|.
$$ 

The local time of SRW is defined by
$$
\boldsymbol{\xi^*}(k,n):=\sum_{i=1}^n I\{{\bf S}_i^*=k\},
\quad k=0,\pm 1,\pm 2,\ldots, \,\, n=1,2,\ldots,
$$
and that of an SBM is defined by
\begg
\boldsymbol{\eta^*}(x,t):=\lim_{\varepsilon\to 0} \frac{1}{2\varepsilon}
\int_0^t I\{{\bf B^*}(s)\in (x-\varepsilon,x+\varepsilon)\}\, ds. 
\label{skewloc}
\endd

Various results for local times are given in Walsh \cite{WAL}, Burdzy and 
Chen \cite{BC}, Lyulko \cite{Ly}, Gairat and Shcherbakov \cite{GS}.  

Here  we  first give a joint strong invariance principle for SBM and SRW
and their local times. 

\begin{thm} A probability space with an {\rm SBM}  
$\{{\bf B^*}(t),\, t\geq 0\}$ 
and an {\rm SRW}  $\{{\bf S}_n^*,\, n=0,1,2,\ldots\}$ on it with the same 
parameter $p$, can be so constructed that, as $n\to\infty$, we have
\begin{equation}
|{\bf S}_n^*-{\bf B^*}(n)|=O(n^{1/4}(\log n)^{1/2}(\log\log n)^{1/4}) \quad a.s.
\label{inv1}
\end{equation}
and 
\begin{equation}
\max_{1\leq k\leq n}\sup_{x\in {\bz}}
|\boldsymbol{\xi^*}(x,k)-\boldsymbol{\eta^*}(x,k)|
=O(n^{1/4}(\log n)^{1/2}(\log\log n)^{1/4})\quad a.s. 
\label{inv2}
\end{equation}
\end{thm}

\medskip\noindent
{\bf Proof.} The statement of (\ref{inv1}) is a special case of Theorem 1.1 
when $N=2$, with $p_1=p$, $p_2=1-p=q$. 

The conclusion of (\ref{inv2}) was established in Cs\"org\H o and Horv\'ath 
\cite{CsHo} in the case of a simple symmetric random walk, i.e., when $p=1/2$. 
We are to conclude (\ref{inv2}) in the present context quite similarly to  
the proof of Theorem 1 in \cite{CsHo} with suitable modifications.

Let $\{{\bb^*}(t),\, t\geq 0\}$ be an SBM and define $\tau_0=0$,
$$\tau_1=\inf\{t:\, t>0,\, |{\bb}^*(t)|=1\},  \quad 
\tau_n=\inf\{t:\, t>\tau_{n-1},\, |{\bb^*}(t)-{\bb^*}(\tau_{n-1})|=1\},$$
$n=2,3,\ldots$. Using these Skorokhod stopping times, in view of (\ref{sdef}) 
and the proof of Theorem 1.1 in \cite{CCFR83}, it is easy to see that 
${\bs}_0^*=0,$ $${\bs}_n^*=\sum_{i=1}^n X_i^*,\, n=1,2\ldots, \,\,{\rm with}\,\,
X_i^*={\bb^*}(\tau_i)-{\bb^*}(\tau_{i-1}),\, i=1,2,\ldots,$$ is an SRW. We note
in passing that, along these lines, we can again conclude the statement
of (\ref{inv1}).

Proving now (\ref{inv2}), put
 $$a_i(x)=
\boldsymbol{\eta^*}(x, \tau_{\nu(i)+1})-\boldsymbol{\eta^*}(x, \tau_{\nu(i)}),
$$
where
$$\nu(1)=\min\{k:\, k\geq 0,\, {\bb^*}(\tau_k)={\bs}_k^*=x\},\, 
\nu(n)=\min\{k:\, k>\nu(n-1),\, {\bb^*}(\tau_k)={\bs}_k^*=x\}, \,\, n\geq 2.$$
Observe that $a_i(x)$,  $i=1,2 \dots$ for any fixed $x\in {\bz}$, are i.i.d. 
random variables. Moreover, even though they are defined using an SBM, 
their common distribution is identical with those, which  one gets  using  
a standard Brownian motion. Consequently, we can use Borodin and Salminen 
\cite{BS}, p. 173, (3.3.2), to conclude that the random variables 
$\{a_i(x),\, i\geq 1\},$ for all $x\in {\bz},$  are i.i.d. 
with exponential density function with parameter $1$. Thus, for further use, 
for $i=1,2, \dots, $and $x\in{\bz}$, we have
\begin{equation}
{\bf P}(a_i(x)\geq y)=e^{-y},\quad y\geq 0.
\label{nahat}
\end{equation}

Observe that, for all $n\geq 1$, we have
\begin{equation}
{\bf P}\left(\left|\sum_{i=1}^{\boldsymbol{\xi^*}(x,n)}a_i(x)-
\boldsymbol{\eta^*}(x,\tau_n)\right|\leq
\boldsymbol{\eta^*}(0,\tau_1)+a_{\boldsymbol{\xi^*}(x,n)}(x)\ 
{\rm for\, \, all} \, x\in {\bz}\right)=1.
\label{sumai}
\end{equation}
Based on (\ref{nahat}), we have 
$$
{\bf P}\left(\max_{1\leq i\leq n}\, \max_{-n\leq x\leq n}a_i(x)>5\log 
n\right)\leq 3n^2 {\bf P}(a_1(1)>5\log n)\leq \frac{1}{n^2}
$$ 
and, therefore, as $n\to\infty,$ by the Borel-Cantelli Lemma, 
\begin{equation}
\max_{1\leq i\leq n}\, \max_{-n\leq x\leq n}a_i(x)=O(\log n)\quad a.s.
\label{maxai}
\end{equation}

\noindent
In view of (\ref{sumai}) and (\ref{maxai}), and observing that 
$\boldsymbol{\eta^*}(0,\tau_1)$ has the same distribution as $a_i(x),$ 
we conclude 

\begin{equation}
\max_{-n\leq x\leq n}\left|\sum_{i=1}
^{\boldsymbol{\xi^*}(x,n)}a_i(x)-\boldsymbol{\eta^*}(x,\tau_n)\right|
=O(\log n)\quad a.s.,
\label{sumai2}
\end{equation}
as $n\to\infty$.

Observe that by the law of the iterated logarithm (LIL) of Kesten
\cite{Ke} for the maximum in $x\in {\bz}$ of the local time of the simple 
symmetric random walk $S(\cdot)$, we have 
\begin{equation}
\limsup_{n\to\infty}(2n\log\log n)^{-1/2}\max_{x\in {\bz}}\xi(x,n)=1\quad a.s.
\label{limsup}
\end{equation}
Consequently, for the maximum local time of ${\bf S}_n^*$, we have
\begin{equation}
\limsup_{n\to\infty}(2n\log\log n)^{-1/2}\max_{x\in {\bz}}
\boldsymbol{\xi^*}(x,n)\leq 2\quad a.s.
\label{limsup2}
\end{equation}

\noindent
Using (\ref{nahat}) we get 
$$ {\bf E}(e^{t(a_i(x)-1)})=\frac{e^{-t}}{1-t}, \qquad   
{\bf E}(e^{t(1-a_i(x))})=\frac{e^{t}}{1+t},$$
so for $0\leq t\leq 1/2$ we have
$$ {\bf E}(e^{t(a_i(x)-1)}) \leq e^{t^2}, \qquad   
{\bf E}(e^{t(1-a_i(x))})\leq e^{t^2}.$$ 
This implies that we can apply the inequality of Petrov \cite{PET} (Section 
III.5, Supplement 10, page 58) for the maximum of partial sums to get that for 
$0\leq z\leq K$
$$
{\bf P}\left(\max_{1\leq k\leq K}
\left|\sum_{i=1}^k(a_i(x)-1)\right| \geq z \right)
$$
$$
\leq {\bf P}\left(\max_{1\leq k\leq K}
\sum_{i=1}^k(a_i(x)-1)\geq z \right)
+{\bf P}\left(\max_{1\leq k\leq K}
\sum_{i=1}^k(1-a_i(x)) \geq z \right)
\leq 2\exp\left(-\frac{z^2}{4K}\right). $$ 
Putting $K=4(n\log\log n)^{1/2}$, we obtain
$$
{\bf P}\left(\max_{1\leq k\leq 4(n\log\log n)^{1/2}}
\left|\sum_{i=1}^k(a_i(x)-1)\right|
>4\sqrt{3}(n\log\log n)^{1/4}(\log n)^{1/2}\right)\leq\frac{2}{n^3},
$$
and hence
$$
{\bf P}\left(\max_{|x|\leq n}\max_{1\leq k\leq 4(n\log\log n)^{1/2}}
\left|\sum_{i=1}^k(a_i(x)-1)\right|>
4\sqrt{3}(n\log\log n)^{1/4}(\log n)^{1/2}\right)\leq \frac{C}{n^2}
$$
with some constant $C>0$. 
Using (\ref{limsup2}), by the Borel-Cantelli Lemma  we obtain
$$
\max_{|x|\leq n}\left|\sum_{i=1}^{\boldsymbol{\xi^*}(x,n)}(a_i(x)-1)\right|
\leq\max_{|x|\leq n}\max_{1\leq k\leq 4(n\log\log n)^{1/2}}
\left|\sum_{i=1}^k(a_i(x)-1)\right|=O((n\log\log n)^{1/4}(\log n)^{1/2})
$$
almost surely, as $n\to\infty$. So we arrived to the conclusion that
\begin{equation}
\limsup_{n\to\infty}n^{-1/4}(\log\log n)^{-1/4}(\log n)^{-1//2}
\max_{-n\leq x\leq n}\left|\sum_{i=1}^{\boldsymbol{\xi^*}(x,n)}a_i(x)
-\boldsymbol{\xi^*}(x,n)\right|\leq C_1\quad a.s.
\label{maxai2}
\end{equation}
with some positive constant $C_1$.

It is well known that, in terms of the Skorokhod stopping times 
$\{\tau_i,\, i\geq 0\}$, the random variables 
$\{\tau_i-\tau_{i-1},\, i\geq 1\}$ are i.i.d. random variables with 
mean 1 and variance 1. Thus, by the LIL for partial sums, we have
$$
\limsup_{n\to\infty}(2n\log\log n)^{-1/2}|\tau_n-n|=1\quad a.s.
$$ 
The latter LIL combined with Theorem 3 of Cs\'aki et al. \cite{CCFR83}
with $g(n)=4(n\log\log n)^{1/2}$ yields
\begin{equation}
\limsup_{n\to\infty}n^{-1/4}(\log\log n)^{-1/4}(\log n)^{-1/2}
\sup_x|\boldsymbol{\eta^*}(x,n\pm g(n))-\boldsymbol{\eta^*}(x,n)|
\leq C_2\quad a.s.
\label{incr}
\end{equation} 
with some positive constant $C_2$.
 
Also,
\begin{equation}
\boldsymbol{\xi^*}(x,n)=0\quad {\rm if }\, |x|>n,\quad {\rm and}\,\,\,
\lim_{n\to\infty}\sup_{|x|>n}\boldsymbol{\eta^*}(x,n)=0\quad a.s.
\end{equation}

As a consequence of (\ref{sumai2}), (\ref{limsup2})--(\ref{incr}), we now 
conclude
\begin{equation}
\limsup_{n\to\infty}n^{-1/4}(\log\log n)^{-1/4}(\log n)^{-1/2}
\sup_{x\in {\bz}}|\boldsymbol{\xi^*}(x,n)-\boldsymbol{\eta^*}(x,n)|
\leq C_1+C_2\quad a.s.
\end{equation}
that, in turn, implies (\ref{inv2}), and hence also concludes the proof
of Theorem 2.1. \, \, \, $\Box$

\smallskip
\noindent
{\bf Remark}. In \cite{CsHo} Cs\"org\H o and Horv\'ath show that, using 
the Skorokhod embedding scheme, the rate of convergence in their Theorem 1
is best possible. The latter result is based on their conclusion that, 
using the Skorokhod embedding in case of the simple symmetric random walk, 
one has
$$
\limsup_{n\to\infty}n^{-1/4}(\log\log n)^{-1/4}(\log n)^{-1/2}
\max_{1\leq k\leq n}|\xi(0,k)-\eta(0,k)|=C\quad a.s.
$$ 
with some positive constant $C$. Since the local time of an SRW ${\bs}_n^*$ at 
zero and that of an SBM ${\bb}^*(t)$ respectively coincide with those of  
a simple symmetric random walk and standard Brownian motion, the above 
conclusion also implies that the rate of approximation in (\ref{inv2}) of 
our Theorem 2.1 is also best possible when using the Skorokhod embedding 
scheme.  
   
This invariance principle is suitable for proving so-called first order
limit theorems for local times. For example, in view of having the LIL for 
the maximal local time of an SRW as in (\ref{limsup2}), we can conclude the 
same LIL for the maximal local time of an SBM. Establishing the exact 
constant for any one of these two LIL's, the same constant would be inherited 
by the other one in hand.  

For the so-called second order limit theorems we introduce the following 
iterated Brownian motion, or iterated Wiener process. Let $W(t),\, t\geq 0$
be a standard Wiener process on the line, and let $\tilde\eta(t)$ a Wiener 
local time at zero, independent of $W(\cdot)$. Put
$$
Z(t)=W(\tilde\eta(t)),\quad t\geq 0.
$$     
In the sequel we call $Z(t)$ iterated Brownian motion (IBM). For fixed real 
number $x\neq 0$ we define $c(\cdot ),$ for throughout  use from now on, as

\[c(x):=\left \{ \begin{array}{ll} 2p,
 &\mbox
{ $x>0,$}
 \\
 2q,
 &\mbox
{ $x< 0 .$}
 \\
\end{array}
\right. \]

The next strong invariance principles deal with second order limit theorems 
for the local times in hand, at fixed location for large times.

\begin{thm} For any fixed integer $x\neq 0$, a probability space with an 
{\rm SRW} $\{{\bs}^*_n,\, n=0,1,\ldots\}$, and an {\rm IBM} $\{Z_1(t),
\, t\geq 0\}$ on it can be so constructed  that for the local times of 
{\rm SRW}, as $n\to\infty$, we have 
\begin{equation}
\boldsymbol{\xi}^*(x,n)-c(x)\boldsymbol{\xi}^*(0,n)=
(c(x)(4|x|-1)-c^2(x))^{1/2}\,Z_1(n)+O(n^{\mu})\quad a.s. 
\label{inva1}
\end{equation}
with some $7/16<\mu<1/4$.
\end{thm}
\begin{thm} For any fixed real number $x\neq 0$, a probability space with 
an {\rm SBM} $\{{\bb}^*(t),\quad t\geq 0\}$, and an {\rm IBM} 
$\{Z_2(t),\, t\geq 0\}$, can be so constructed that for the local times of 
{\rm SBM}, as $t\to\infty$, we have
$$
\boldsymbol{\eta}^*(x,t)-c(x)\boldsymbol{\eta}^*(0,t)=
2(c(x)|x|)^{1/2}Z_2(t)+O(t^{\mu}) \quad a.s.
$$
with some $7/16<\mu<1/4$.
\end{thm} 

\noindent{\bf Proof of Theorem 2.2.} Define
$$
\rho_0=0,\, \rho_k=\min\{i>\rho_{k-1}:\, \bs^*_i=0\},\, k=1,2,\ldots,
$$
the return times of SRW to the origin. In fact, $\rho_k,\, k=1,2,\ldots,$
are also the return times of the simple random walk $S(n)$ that was used 
to construct SRW in (2.1). Then, for fixed integer $x\neq 0$ we have 
$$
\bp(\boldsymbol{\xi}^*(x,\rho_1)=0)=1-\frac{c(x)}{2|x|},
$$
$$
\bp(\boldsymbol{\xi}^*(x,\rho_1)=m)=
\frac{c(x)}{4x^2}\left(\frac{2|x|-1}{2|x|}\right)^{m-1},
\quad m=1,2,\ldots
$$
The two formulas above are simple consequences of Theorem 9.7 of 
R\'ev\'esz  (\cite{R13}, page 102), where the above distributions are given 
for a simple symmetric random walk, i.e., for $p=1/2$. Our conclusions are 
modified, as the skew random walk has unequal probabilities for stepping 
right or left from the origin, yielding $c(x)$ in our formulas. 
Accordingly, one can obtain
$$
{\bf E}(\boldsymbol{\xi}^*(x,\rho_1))=c(x),\quad 
Var\boldsymbol{\xi}^*(x,\rho_1)=c(x)(4|x|-1)-c^2(x).
$$
Put
$$
V_i=\boldsymbol{\xi}^*(x,\rho_i)-\boldsymbol{\xi}^*(x,\rho_{i-1}),
\quad i=1,2,\ldots
$$
Then $V_i$ are i.i.d. random variables with distributions given above
with finite moment generating function. So by strong invariance results of 
Koml\'os et al. \cite{KMT} we may conclude
\begin{equation}
\boldsymbol{\xi}^*(x,\rho_r)-c(x)r=(c(x)(4|x|-1)-c^2(x))^{1/2}W(r)+O(\log r) 
\quad a.s.
\label{KMT}
\end{equation}
as $r\to\infty$, where $W(\cdot)$ is a standard Wiener process. On noting that,
by the definition of Spitzer \cite{SP} (cf. also Cs\'aki et al. 
\cite{CCFR92}), the simple random walk $S(n)$ is aperiodic, from here on we 
can follow the method and proof of \cite{CCFR92} to prove Theorem 2.2. We can 
exactly proceed as in \cite{CCFR92}, since in our case of skew random walk, 
the local time $\boldsymbol\xi^*(x,\rho_n)$ can be treated as the sum of i.i.d.
random variables $V_i$, where the excursions $(\rho_{i-1},\rho_i)$,
$i=1,2,\ldots,$ are in fact those of a simple random walk. By Lemma 2.2 in
\cite{CCFR92} we can construct a positive stable 1/2 process (inverse 
Brownian local time) $U_n$ such that
$$
|\rho_n-U_n|=O(n^{5/3})\quad a.s.
$$  
as $n\to\infty$. But constructing  $W(\cdot)$ and $U_{\cdot}$ jointly, 
as above, they are not independent. Following the method and proof of 
Proposition 1 in \cite{CCFR92} we can construct independent standard Wiener 
process $W^{(2)}(\cdot)$ and inverse Brownian local time process 
$U^{(1)}_{\cdot}$ so that 
$$
\boldsymbol{\xi}^*(x,\rho_n)-c(x)n=
(c(x)(4|x|-1)-c^2(x))^{1/2}W^{(2)}(n)+O(n^{2\psi})\quad a.s.
$$
and
$$
|\rho_n-U^{(1)}_n|=O(n^{4\mu}) \quad a.s.,
$$
as $n\to\infty$, where $\mu<\psi/2$, and 
$$
\frac{1}{8}<\psi<\frac{1}{4},\qquad \frac{7}{16}<\mu<\frac{1}{2}.
$$
Note that $\delta$ in (2.8) and (2.9) of \cite{CCFR92} can be choosen 
arbitrary large, since $\boldsymbol{\xi}^*(x,\rho_1)$ have moment generating 
function. In particular, $\mu=\psi=15/62$ can be choosen. Now, just as Theorem 
1 follows from Proposition 1 in \cite{CCFR92}, we can similarly obtain (\ref{inva1}) 
in Theorem 2.2. \, \, \, \, $\Box$   
 
\medskip
\noindent{\bf Proof of Theorem 2.3.}
To prove Theorem 2.3, consider the skew Brownian motion ${\bb}^*(\cdot)$ as a 
diffusion process with scale function  
$$
s(x)=\frac{x}{c(x)},\quad x\neq 0,
$$
and speed measure 
$$
m(dx)=2c(x)dx,\quad x\neq 0.
$$
We note that in It\^o and McKean \cite{IM}, and also in Borodin and Salminen 
\cite{BS}, Appendix 1.12, the scale function is given as twice of the above 
scale function and the speed measure as half of the above speed measure, but 
this is equivalent to the quantities given above, in that they yield the usual 
scale and speed for Brownian motion, i.e., for $p=1/2$. 

In order to apply a  result of Cs\'aki and Salminen \cite{CsS}, we note that 
the local time $L_t^x$ is defined there with respect to the speed measure, 
while $\boldsymbol{\eta}^*(x,t)$ is defined in (\ref{skewloc}) with respect 
to the Lebesgue measure. One can easily see that
$$
L_t^x=\frac{1}{2c(x)}\boldsymbol{\eta}^*(x,t),\quad x\neq 0,
$$
$$
L_t^0=\frac{1}{2}\boldsymbol{\eta}^*(0,t). 
$$
Let $\tau:=\min\{s:\, \boldsymbol{\eta}^*(0,s)=1\}$. Using the formula (2.18) 
in \cite{CsS} for $t=1/2$, we obtain for $\beta>0$,
$$
{\bf E}\left(\exp\left(-\beta\boldsymbol{\eta}^*(x,\tau)\right)\right)=
\exp\left(-\frac{c(x)\beta}{1+2\beta |x|}\right), \quad x\neq 0,
$$
from which
$$
{\bf E}(\boldsymbol{\eta}^*(x,\tau))=c(x),\quad Var(\boldsymbol{\eta}^*(x,\tau))
=4c(x)|x|, \quad x\neq 0.
$$
 Theorem 2.3 can now be proved similarly to Theorem 2.2, referring to 
\cite{CCFR92}. Let $B(\cdot)$ be the standard Brownian motion that is used in 
construction (\ref{sdef2})  in \cite{CCFR92}, and let $U(\cdot)$ be its 
inverse local time at zero. Put
$$
Y_i=\boldsymbol{\eta}^*(x,U(i))-\boldsymbol{\eta}^*(x,U(i-1)),\, \, 
i=1,2,\ldots
$$
Then
$$
\boldsymbol{\eta}^*(x,U(n))=\sum_{i=1}^n Y_i,
$$
which is a sum of i.i.d. random variables that can be approximated by 
a Wiener process as in (\ref{KMT}). Now using the method and proof of Theorem 2 
and Proposition 2 in \cite{CCFR92}, we obtain Theorem 2.3 as asserted  above. 
\, \, \, \, $\Box$

Next we state some consequences of the strong approximations given in 
Theorems 2.2 and 2.3. First, we coclude the following Dobrushin
\cite{Do}-type theorems, where ${\buildrel{d}\over\longrightarrow}$ means 
convergence in distribution.
\begin{corollary}
For fixed integer $x\neq 0$ we have, as $n\to\infty$, 
$$
\frac{\boldsymbol{\xi}^*(x,n)-c(x)\boldsymbol{\xi}^*(0,n)}
{(c(x)(4|x|-1)-c^2(x))^{1/2}n^{1/4}}\, {\buildrel{d}\over\longrightarrow}\,
 U\sqrt{|V|}, 
$$
where $U$ and $V$ are independent standard normal random variables.
\end{corollary}
\begin{corollary}
For any fixed real number $x\neq 0$ we have, as $t\to\infty$,
$$
\frac{\boldsymbol{\eta}^*(x,t)-c(x)\boldsymbol{\eta}^*(0,t)}
{(4c(x)|x|)^{1/2}t^{1/4}}\, {\buildrel{d}\over\longrightarrow}\, 
U\sqrt{|V|}, 
$$
where $U$ and $V$ are independent standard normal random variables.
\end{corollary}

Strassen-type theorems for the IBM, that also imply limsup results, are 
given in Cs\'aki et al. \cite{CCFR95}, \cite{CFR}, and Hu et al. \cite{HPS}. 
A Chung-type liminf result is shown in Hu et al. \cite{HPS}. Applying these
results, we also have the following corollaries to Theorems 2.2 and  2.3, 
respectively.

\begin{corollary}
For  any fixed integer $x\neq 0,$ the set of functions
$$
g_n(s)=\frac{\boldsymbol{\xi}^*(x,[sn])-c(x)\boldsymbol{\xi}^*(0,[sn])}
{(c(x)(4|x|-1)-c^2(x))^{1/2}2^{5/4}3^{-3/4}n^{1/4}(\log\log n)^{3/4}}, \quad
0\leq s\leq 1,
$$
is relatively compact in $C[0,1]$ and the set of its limit points,
as $n\to \infty,$ is the set of functions $f(s),\, 0\leq s\leq 1$, with
$f(0)=0$, that are absolutely continuous with respect to Lebesgue measure
and  
$$
\int_0^1|\dot f(s)|^{3/4}\, ds\leq 1.
$$
Consequently,
$$
\limsup_{n\to\infty}\frac{\boldsymbol{\xi}^*(x,n)-c(x)\boldsymbol{\xi}^*(0,n)}
{(c(x)(4|x|-1)-c^2(x))^{1/2}2^{5/4}3^{-3/4}n^{1/4}(\log\log n)^{3/4}}=1
$$
almost surely. Also,
$$
\liminf_{n\to\infty}\frac{(\log\log n)^{3/4}}
{(c(x)(4|x|-1)-c^2(x))^{1/2}n^{1/4}}
\max_{0\leq k\leq n}|\boldsymbol{\xi}^*(x,k)-c(x)\boldsymbol{\xi}^*(0,k)|
=\left(\frac{3\pi^2}{8}\right)^{3/4}
$$
almost surely.
\end{corollary}

\begin{corollary}
For any fixed real number $x\neq 0$ the set of functions
$$
g_t(s)=\frac{\boldsymbol{\eta}^*(x,st)-c(x)\boldsymbol{\eta}^*(0,st)}
{2(c(x)|x|)^{1/2}2^{5/4}3^{-3/4}t^{1/4}(\log\log t)^{3/4}}, \quad
0\leq s\leq 1,
$$
is relatively compact in $C[0,1]$ and the set of its limit points,
as $t\to \infty,$ is the set of functions $f(s),\, 0\leq s\leq 1$, with
$f(0)=0$, that are absolutely continuous with respect to Lebesgue measure
and  
$$
\int_0^1|\dot f(s)|^{3/4}\, ds\leq 1.
$$
Consequently,
$$
\limsup_{t\to\infty}\frac{\boldsymbol{\eta}^*(x,t)-c(x)\boldsymbol{\eta}^*(0,t)}
{2(c(x)|x|)^{1/2}2^{5/4}3^{-3/4}t^{1/4}(\log\log t)^{3/4}}=1
$$
almost surely. Also,
$$
\liminf_{t\to\infty}\frac{(\log\log t)^{3/4}}
{2(c(x)|x|)^{1/2}t^{1/4}}
\max_{0\leq s\leq t}\boldsymbol|{\eta}^*(x,s)-c(x)\boldsymbol{\eta}^*(0,s)|
=\left(\frac{3\pi^2}{8}\right)^{3/4}
$$
almost surely.
\end{corollary}
  
\subsection{Local time on spider}
First we give a joint invariance principle for the  Brownian spider 
and the random walk on spider ${\bf SP}(N)$, and their local times, an analogue 
of Theorem 2.1.

\begin{thm} On a rich enough probability space one can define a Brownian
spider $\{{\bf B}(t),\, t\geq 0\}$ and random walk on the spider
$\{{\bf S}_n,\, n=0,1,2,\ldots\}$, both on ${\bf SP}(N)$, and both
selecting their legs  with the same  probabilities
$p_j,\, j=1,2,\ldots,N,$ so that, as $n\to\infty$, we have
$$
|{\bf S}_n-{\bf B}(n)|=O(n^{1/4}(\log n)^{1/2}(\log\log n)^{1/4}) \quad a.s.
$$
and 
$$
\sup_{(r,j)}|\boldsymbol{\xi}((r,j),n)-\boldsymbol{\eta}((r,j),n)|
=O(n^{1/4}(\log n)^{1/2}(\log\log n)^{1/4})\quad a.s.,
$$
where $\sup$ is taken for $r=0,1,2,\ldots$, $j=1,2,\ldots, N$.
\end{thm}

\noindent{\bf Proof.} The proof can be reduced to the case $N=2$. Namely, 
for a fixed $j,$ consider the SRW and SBM with $p=p_j$, $q=1-p_j$. For this $j$,
Theorem 2.4 follows from Theorem 2.1. Repeating the proof for all 
$j=1,\ldots, N$, yields the proof of Theorem 2.4.\ \ \ \ \ $\Box$ 

Concerning the so-called second order limit theorems for the local time,
we have the following results, corresponding respectively to Theorems 2.2 
and 2.3, on replacing $c(x)$ by $2p_j$.  

\begin{thm}
For any fixed positive integer $x$ and  each $j=1,2,\ldots, N$, a probability 
space with an {\rm RWS} ${\bf S}_n $ and an {\rm IBM} $Z_3(t)$ can be so 
constructed that, as $n\to\infty$, we have 
$$
\boldsymbol{\xi}((x,j),n)-2p_j\boldsymbol{\xi}(0,n)=
(2p_j(4x-1)-4p_j^2)^{1/2}\,Z_3(n)+o(n^{\mu}) \quad a.s.,
$$
where $\mu$ is as in Theorem {\rm 2.2.}
\end{thm}

\begin{thm}
For any fixed positive real number $x$ and  each $j=1,2,\ldots, N,$ a 
probability space with an {\rm BMS} ${\bf B}(t)$ and an {\rm IBM} 
$Z_4(t)$ can be so constructed that, as $t\to\infty$, we have 
$$
\boldsymbol{\eta}((x,j),t)-2p_j\boldsymbol{\eta}(0,t)=
(8p_jx)^{1/2}\,Z_4(t)+o(t^{\mu}) \quad a.s.,
$$
where $\mu$ is as in Theorem {\rm 2.3.}
\end{thm}

Similarly to Corollaries 2.1-2.4, we have the following Dobrushin and 
Strassen type respective consequences of Theorems 2.5 and 2.6 for local
times on a spider. 
\begin{corollary}
For any fixed integer $x>0$ and  each $j=1,2,\ldots, N,$ we have, as 
$n\to\infty$,
$$
\frac{\boldsymbol{\xi}((x,j),n)-2p_j\boldsymbol{\xi}(0,n)}
{(2p_j(4x-1)-4p_j^2)^{1/2}n^{1/4}}\, {\buildrel{d}\over\longrightarrow}\,
 U\sqrt{|V|},
$$
where $U$ and $V$ are independent standard normal random variables.
\end{corollary}
\begin{corollary}
For any fixed real number $x>0$ and  each $j=1,2,\ldots, N,$ we have, 
as $t\to\infty$,
$$
\frac{\boldsymbol{\eta}((x,j),t)-2p_j\boldsymbol{\eta}(0,t)}
{(8p_jx)^{1/2}t^{1/4}}\, {\buildrel{d}\over\longrightarrow}\,
U\sqrt{|V|},
$$
where $U$ and $V$ are independent standard normal random variables.
\end{corollary}
\begin{corollary}
For any fixed integer $x>0$ and  each $j=1,2,\ldots, N,$ the set of functions
$$
g_n(s)=\frac{\boldsymbol{\xi}((x,j),[sn])-2p_j\boldsymbol{\xi}(0,[sn])}
{(2p_j(4x-1)-4p_j^2)^{1/2}2^{5/4}3^{-3/4}n^{1/4}(\log\log n)^{3/4}}, \quad
0\leq s\leq 1,
$$
is relatively compact in $C[0,1]$ and the set of its limit points,
as $n\to \infty,$ is the set of functions $f(s),\, 0\leq s\leq 1$, with
$f(0)=0$, that are absolutely continuous with respect to Lebesgue measure
and
$$
\int_0^1|\dot f(s)|^{3/4}\, ds\leq 1.
$$
Consequently,
$$
\limsup_{n\to\infty}\frac{\boldsymbol{\xi}((x,j),n)-2p_j\boldsymbol{\xi}(0,n)}
{(2p_j(4x-1)-4p_j^2)^{1/2}2^{5/4}3^{-3/4}n^{1/4}(\log\log n)^{3/4}}=1
$$
almost surely. Also,
$$
\liminf_{n\to\infty}\frac{(\log\log n)^{3/4}}
{(2p_j(4x-1)-4p_j^2)^{1/2}n^{1/4}}
\max_{0\leq k\leq n}|\boldsymbol{\xi}((x,j),k)-2p_j\boldsymbol{\xi}(0,k)|
=\left(\frac{3\pi^2}{8}\right)^{3/4}
$$
almost surely.
\end{corollary}
\begin{corollary}
For any fixed real number $x>0$ and  each $j=1,2,\ldots, N,$ the set of 
functions
$$
g_t(s)=\frac{\boldsymbol{\eta}((x,j),st)-2p_j\boldsymbol{\eta}(0,st)}
{(8p_jx)^{1/2}2^{5/4}3^{-3/4}t^{1/4}(\log\log t)^{3/4}}, \quad
0\leq s\leq 1,
$$
is relatively compact in $C[0,1]$ and the set of its limit points,
as $t\to \infty,$ is the set of functions $f(s),\, 0\leq s\leq 1$, with
$f(0)=0$, that are absolutely continuous with respect to Lebesgue measure
and
$$
\int_0^1|\dot f(s)|^{3/4}\, ds\leq 1.
$$
Consequently,
$$
\limsup_{t\to\infty}\frac{\boldsymbol{\eta}((x,j),t)-2p_j\boldsymbol{\eta}(0,t)}
{(8p_jx)^{1/2}2^{5/4}3^{-3/4}t^{1/4}(\log\log t)^{3/4}}=1
$$
almost surely. Also,
$$
\liminf_{t\to\infty}\frac{(\log\log t)^{3/4}}
{(8p_jx)^{1/2}t^{1/4}}
\max_{0\leq s\leq t}|\boldsymbol{\eta}((x,j),s)-2p_j\boldsymbol{\eta}(0,s)|
=\left(\frac{3\pi^2}{8}\right)^{3/4}
$$
almost surely.
\end{corollary}

\section{Occupation times.}
\renewcommand{\thesection}{\arabic{section}} \setcounter{equation}{0}
\setcounter{thm}{0} \setcounter{lemma}{0}

In this section we consider the time the RWS, or  the BMS, spends 
on particular legs. First we give the definitions of the occupation times. 
Recall that $L_j$ denotes the leg $j$ of the spider $\textbf{SP}(N)$.
Let $\bs_j$, $j=0,1,\ldots,$ be a random walk on the spider $\textbf{SP}(N)$,
with probabilities $p_1,\ldots,p_N$. The occupation time on leg $j$ in
$n$ steps is defined by
$$
T(j,n):=\sum_{k=1}^nI\{\bs_k\in L_j\}.
$$   

\noindent
The occupation time of the Brownian spider $\bb(s)$, $s\geq 0,$ up to time 
$t$ is defined by
$$
Z(j,t):=\int_0^t I\{\bb(s)\in L_j\}\, ds.
$$

\begin{thm} On the probability space of {\rm Theorem 2.4}  as 
$n\to\infty$, we have for any $\epsilon>0,$
$$
\sup_{1\leq j\leq N} |T(j,n)-Z(j,n)|=O(n^{3/4+\varepsilon})\quad a.s.
$$
\end{thm}

\noindent{\bf Proof of Theorem 3.1.}
The occupation time of RWS on a particular leg of $\textbf{SP}(N)$ can be 
considered as having only two legs, one to which the walker proceeds with 
probability $p_j$ from the origin, and another one to which it goes with 
probability $1-p_j$. Based on this idea, we consider the following setup.

Let $N=2,$  $p=p_1,$\ and $q=p_2=1-p$.
This corresponds to a skew random walk $\{{\bf S}^*_i,\, i=0,1,\ldots\}$ with
transition probabilities
$$
{\bf P}(0,1)=p, \,\, {\bf P}(0,-1)=q=1-p,\,\, 
{\bf P}(x,x+1)={\bf P}(x,x-1)=1/2,\,\, x\neq 0.
$$

\noindent
Let $T^*(1,n)$ be the occupation time of the positive half-line
of a skew random walk, i.e.,
$$
T^*(1,n)=\sum_{i=1}^n I\{{\bf S}^*_i>0\}
=\sum_{i=1}^\infty\boldsymbol{\xi^*}(i,n).
$$
Moreover, let $Z^*(1,t)$ be the occupation time of the positive half-line
of a skew Brownian motion, i.e.,
$$
Z^*(1,t)=\int_0^t I\{B^*(s)> 0\}ds=
\int_0^\infty\boldsymbol{\eta^*}(x,t)\, dx. 
$$

\medskip
The limiting density of $T^*(1,n)$ can be obtained from the density of 
$Z^*(1,t)$, for which we have
$$
{\bf P}\left(\frac{Z^*(1,t)}{t}\in dx\right)=
\frac{pq}{\pi\sqrt{x(1-x)}(p^2(1-x)+q^2x)}\, dx,\quad 0<x<1.
$$ 
 
The latter formula is given in Lamperti
\cite{LA56} as a limiting distribution of occupation time
of certain discrete time processes, including SRW. For this formula and related 
results, see also Appuhamillage et al. \cite{ABT}.

So we have to prove that for a suitable construction of SBM and SRW, as 
$n\to\infty,$ we have for any $\varepsilon>0$, 
$$
|T^*(1,n)-Z^*(1,n)|=O(n^{3/4+\varepsilon})\quad a.s.,
$$
It follows from the second statement of Theorem 2.1 and the law of the 
iterated logarithm that, as $n\to\infty,$
\begin{equation}
\sum_{i=1}^\infty|\boldsymbol\xi^*(i,n)-\boldsymbol\eta^*(i,n)|
=O(n^{3/4+\varepsilon})\quad a.s. 
\label{invar}
\end{equation}
 Also,
$$
\left|Z^*(1,n)-\sum_{i=1}^\infty\boldsymbol\eta^*(i,n)\right|
\leq \sum_{i=1}^\infty \sup_{i-1<x\leq i}
|\boldsymbol\eta^*(x,n)-\boldsymbol\eta^*(i,n)|.
$$
Now, for Wiener local time Bass and Griffin \cite{BG}, Lemma 5.3, proved that
for every $\varepsilon>0$
$$
\sup_k\, \sup_{s\leq t}\, \sup_{z\in [k,k+1]}
|\eta(z,s)-\eta(k,s)|=o(t^{1/4+\varepsilon/2})
$$
almost surely, as $t\to\infty$. The skew Brownian motion can be obtained from 
standard Brownian motion by reflecting some excursions from negative to 
positive part or from positive to negative part. Hence for $i-1<x\leq i$,
$$
|\boldsymbol\eta^*(x,n)-\boldsymbol\eta^*(i,n)|\leq
|\eta(x,n)-\eta(i,n)|+|\eta(-x,n)-\eta(-i,n)|,
$$ 
and thus, by the law of the iterated logarithm
$$
\left|Z^*(1,n)-\sum_{i=1}^\infty\boldsymbol\eta^*(i,n)\right|
=O(n^{3/4+\varepsilon})\quad a.s.,
$$
as $n\to\infty$. This, combined with (\ref{invar}), proves 
Theorem 3.1. \, \, \, $\Box$

Since the occupation time is usually of order $n$, using this strong 
invariance, it suffices to prove strong theorems either for RWS or BMS.  

For the joint distribution of $\{Z(j,t)/t,\, j=1,\ldots,N\}$, Barlow et al. 
\cite {BPY2} (see also Yano \cite{Ya}) have shown the following equality in 
distribution: up to any time $t>0,$
\begg
\left \{\frac{Z(j,t)}{t}, \quad j=1,\ldots,N\right  \} {\buildrel d \over =} 
\left\{\frac{p_j^2 U_j}{\sum_{k=1}^N p_k^2U_k}
,\quad j=1,\ldots,N\right\},
\label {bl}
\endd
where $U_1,\ldots, U_N$ are independent one-sided stable 1/2 random variables.
This, via Theorem 3.1,  also yields  the joint limiting distribution of 
$\{T(j,n)/n,\, j=1,\ldots,N\},$ as $n\to \infty.$

Now let
$$
T_M(n)=\max_{1\leq j\leq N}T(j,n),\qquad
T_m(n)=\min_{1\leq j\leq N}T(j,n).
$$
and
$$
Z_M(t)=\max_{1\leq j\leq N}Z(j,t),\qquad
Z_m(t)=\min_{1\leq j\leq N}Z(j,t).
$$
For the limsup of $T_M$, $Z_M$ and liminf of $T_m$, $Z_m$, we show that the
Chung-Erd\H os  \cite{CHE52} result for simple symmetric walk remains valid.

\begin{thm}
Let $f(x)$ be a positive nondecreasing function for which
$\lim_{x\to\infty}f(x)=\infty$, $x/f(x)$ is nondecreasing and
$\lim_{x\to\infty}x/f(x)=\infty$. 

Let 
$$
I(f):=\int_1^\infty \frac{dx}{x(f(x))^{1/2}}. 
$$
Then

$$
{\bf P}\left(T_M(n)>n\left(1-\frac{1}{f(n)}\right) i.o. \,\,as \,\, n\to  
\infty \right)=0 \,\,{\rm or}\,\,1,
$$

$$
{\bf P}\left(Z_M(t)>t \left(1-\frac{1}{f(t)}\right) i.o. \,\,as \,\, t\to  
\infty \right)=0 \,\,{\rm or}\,\,1,
$$

$$
{\bf P}\left(T_m(n)<\frac{n}{f(n)}\,\,i.o. \,\,as \,\, n\to  \infty \right)=0 
\,\,{\rm or}\,\,1
$$
and
$$
{\bf P}\left(Z_m(t)<\frac{t}{f(t)}\,\,i.o. \,\,as \,\, t\to  \infty \right)=0 
\,\,{\rm or}\,\,1
$$
according as $I(f)$ converges or diverges.

\end{thm}

\noindent
{\bf Proof of Theorem 3.2}. We will only prove the third  statement, the 
proofs for the others go similarly.  The limiting distribution of $T(i,n)$
for fixed $i$ is the same as that of skew random walk, which is a general
arcsin law (see Lamperti \cite{LA56} or Watanabe, \cite{WAT}). For small
$x$ it is of the same order as that of the arcsin law. Therefore,
$$
c_1\sqrt{x}\leq {\bf P}(T(i,n)\leq nx)\leq c_2\sqrt{x}
$$
with some positive constants $c_1$ and $c_2$. As
\begin{eqnarray*} {\bf P}(T(1,n)\leq nx)\leq {\bf P}(T_m(n)\leq nx)
&=&{\bf P}(\bigcup_{i=1}^N \{ T(i,n)\leq nx\})\\
&\leq& N(\max_{i\leq N} {\bf P}(T(i,n)\leq nx),
\end{eqnarray*}
 it follows that for small  $x$ we have
$$
c_3\sqrt{x}\leq {\bf P}(T_m(n)\leq nx)\leq c_4\sqrt{x}
$$
with positive constants $c_3,c_4$. Based on the proof of Chung-Erd\H os  
\cite{CHE52}, we can prove the theorem as follows.

\noindent
{\bf Convergent part}: Let $n_k=2^k.$  Then we know that for a positive 
nondecreasing function $f(x)$, the integral 
$I(f)=\displaystyle{\int_1^\infty \frac{dx}{x(f(x))^{1/2}}}$ converges
if and only if the sum 
$\displaystyle{\sum_{k=1}^{\infty} \frac{1}{(f(2^k))^{1/2}}}$
converges.

\noindent
First we show that if $I(f)<\infty,$ then
\begg\sum_{k=1}^{\infty} {\bf P}\left(T_m(n_{k-1})\leq \frac{n_k}{f(n_k)}\right)
< {\infty}. \label{konv}\endd

We have
\begin{eqnarray*}
{\bf P}\left(T_m(n_{k-1})\leq\frac{2^k}{f(2^k)}\right)
&=&{\bf P}\left(T_m(2^{k-1})\leq 2\frac{2^{k-1}}{f(2^{k-1})}
\frac{f(2^{k-1})}{f(2^k)}\right)\\
&\leq&{\bf P}\left(T_m(2^{k-1})\leq 2\frac{2^{k-1}}{f(2^{k-1})}\right)
\leq c_2 \sqrt{2}\frac{1}{(f(2^{k-1}))^{1/2}},
\end{eqnarray*}
 where we used that $f(\cdot)$ is nondecreasing. By assumption, we now have
that (\ref{konv}) is convergent. Thus,  for $k$ big enough,
 $$T_m(n_{k-1})\geq \frac{n_k}{f(n_k)}\quad  a.s.$$
 Now for $n_{k-1}\leq n\leq n_k, $  we have for $n$ big enough that
 $$T_m(n)\geq T_m(n_{k-1}) \geq \frac{n_k}{f(n_k)}\geq \frac{n}{f(n)}. $$

\noindent
{\bf Divergent part}: First we show that the theorem is valid for a skew
random walk. To this end  suppose that $N=2,$ and $p_1=p\leq 1/2$ and
$p_2=q=1-p.$ Clearly if $p=1/2$ there is nothing to prove, in view of the
original Chung-Erd\H{o}s theorem. That is to say, we know for the time
$T(1,n)$ spent on the positive side (spent on the first leg), we have
$\displaystyle{{\bf P}(T(1,n)<n/f(n)\,\,i.o.)=1},$
 whenever $I(f)=\infty.$
 If $p<1/2$, then we start with a simple symmetric walk again, and 
keep each excursion spent on the positive side with probability $2p$ and,
with probability $1-2p$, we flip it to the negative side. Thus we get a skew
randow walk with $p_1=p<1/2.$ Denote the time spent on the positive side of 
this walk by $T^p(1,n).$ By construction,
$T^p(1,n)\leq T(1,n),$  so $\displaystyle{{\bf P}(T^p(1,n)<n/f(n)\,\,i.o.)=1},$
whenever $I(f)=\infty.$ Similarly, if we have a skew random walk
with  $p_1=p>1/2, $ we start again with a simple symmetric walk, but now
apply the Chung-Erd\H{o}s theorem for the time spent on the negative side 
(call it leg two). Then we have 
$\displaystyle{{\bf P}(T(2,n)<n/f(n)\,\,i.o.)=1},$whenever $I(f)=\infty.$ 
Now construct the new skew random walk by
keeping each excursion on leg two with probability $2q$, and flipping them
to the first leg with probability $1-2q.$ As before, $T^p(2,n)\leq T(2,n),$
hence $\displaystyle{{\bf P}(T^p(2,n)<n/f(n)\,\,i.o.)=1},$ whenever 
$I(f)=\infty.$ Consequently, for a skew random walk $T^p_m(n)=
\min\{T^p(1,n),T^p(2,n)\}$, we have
$\displaystyle{{\bf P}(T^p_m(n)<n/f(n)\,\,i.o.)=1},$ whenever 
$I(f)=\infty.$ Now for the spider random walk, it is enough to observe 
that selecting the leg (or one of the legs) with 
$p*:=\min_{1\leq i\leq N}\,\, p_i$ and 
putting all the excursions from the remaining legs to a second leg, we create  
a skew random walk for which 
$\displaystyle{{\bf P}(T_m^{p}(1,n)<n/f(n)\,\,i.o.)=1},$
whenever $I(f)=\infty.$ Clearly $T_m(n)\leq T^{p*}(1,n)$, so 
the theorem is now proved. $\Box.$

\medskip
Recall the respective definitions of $T_M$, $Z_M$ and $T_m$, $Z_m$, given 
right before Theorem 3.2.
For the liminf of $T_M$, $Z_M$, and the limsup of $T_m$, $Z_m$, we have
the following results.
\begin{thm}
\begin{equation}
\liminf_{n\to\infty}\frac{T_M(n)}{n}=
\limsup_{n\to\infty}\frac{T_m(n)}{n}=\frac{1}{N}\quad a.s. \label{first}
\end{equation}
and, similarly,
\begin{equation}
\liminf_{t\to\infty}\frac{Z_M(t)}{t}=
\limsup_{t\to\infty}\frac{Z_m(t)}{t}=\frac{1}{N}\quad  a.s.\label{second}
\end{equation}
\end{thm}
\medskip\noindent
{\bf Proof of Theorem 3.3}. We show this result for $Z_m$ and $Z_M$. 
The strong invariance result of Theorem 3.1 implies the conclusion for 
$T_m$ and $T_M$. Recalling (\ref{bl}), we know that therein,  
$U_j, j=1,2,...,N,$ are independent stable 1/2 random variables. Thus  
the events 
$$
\left\{\frac{1-\varepsilon}N\leq \min_{1\leq j\leq N}
\frac{p_j^2 U_j}{\sum_{k=1}^N p_k^2 U_k}\leq \frac{1}{N}\right\}
$$
and
$$
\left\{\frac{1}N\leq \max_{1\leq j\leq N}
\frac{p_j^2 U_j}{\sum_{k=1}^N p_k^2 U_k}\leq \frac{1+\varepsilon}{N}\right\}
$$
have positive probability for $0<\varepsilon<1$. 

Hence,via (3.2), for some $\alpha>0$ and all $t>0, $ we have
$$
{\bf P}\left(\frac{1}{N}\leq\frac{Z_M(t)}{t}\leq 
\frac{1+\varepsilon}{N}\right)\geq \alpha 
$$
and
$$
{\bf P}\left(\frac{1-\varepsilon}{N}\leq\frac{Z_m(t)}{t}\leq 
\frac{1}{N}\right)\geq \alpha. 
$$
Consequently, the events
\begg
\frac{1}{N}
\leq\liminf_{t\to\infty}\frac{Z_M(t)}{t}\leq
\frac{1+\varepsilon}{N} 
\label{haho1} 
\endd
and
\begg
\frac{1-\varepsilon}{N}
\leq\limsup_{t\to\infty}\frac{Z_m(t)}{t}\leq
\frac{1}{N} 
\label{haho2} 
\endd
have positive probability, for all $0<\varepsilon<1$. It follows from the 
0-1 law for Brownian spider (cf. Evans \cite{EV}, Theorem 1.8) that the 
respective events as in (\ref{haho1}) and (\ref{haho2}) hold true  
almost surely. Since $\varepsilon$ is arbitrary, (\ref{second}) follows. 
\ \ \ \ $\Box$

\smallskip\noindent
{\bf Acknowledgements} We sincerely wish to thank the referee of our 
submission for careful reading our manuscript, and for making a number 
of insightful comments and suggestions that helped and prompted us to
improve the presentation and proofs of our results when revising this paper 
for publication.

\end{document}